\documentclass{birkjour}

\usepackage{amsmath,amsthm,amsfonts,amscd,eucal}
\usepackage{xypic}
\usepackage{graphicx}
\usepackage{amssymb}
\usepackage{hyperref}
\usepackage{soul}
\usepackage{microtype}

\def\cb{{\mathcal B}}

\def\ce{{\mathcal E}}
\def\cf{{\mathcal F}}

\def\ch{{\mathcal H}}

\def\ck{{\mathcal K}}

\def\cs{{\mathcal S}}

\def\ga{{\mathfrak A}} 
\def\gb{{\mathfrak B}}\def\gpb{{\mathfrak b}}

\def\gam{{\mathfrak M}}

\def\bc{{\mathbb C}}
\def\bd{{\mathbb D}}
\def\be{{\mathbb E}}
\def\bi{{\mathbb I}}
\def\bl{{\mathbb L}}
\def\bj{{\mathbb J}}
\def\bn{{\mathbb N}}
\def\bp{{\mathbb P}}
\def\bz{{\mathbb Z}}

\def\a{\alpha}

\def\g{\gamma}  \def\G{\Gamma}
\def\d{\delta}  

\def\eps{\varepsilon}

\def\l{\lambda} \def\L{\Lambda}

\def\s{\sigma} 
\def\t{\tau}
\def\f{\varphi}  \def\F{\Phi}
\def\th{\theta}  
\def\om{\omega} \def\Om{\Omega}
\def\z{\zeta}

 \newtheorem{thm}{Theorem}[section]
 \newtheorem{cor}[thm]{Corollary}
 
 \newtheorem{prop}[thm]{Proposition}
 \theoremstyle{definition}
 \newtheorem{defn}[thm]{Definition}
 \theoremstyle{remark}

 \numberwithin{equation}{section}

\begin{document}

%
%
%
%
%
\title[Non-commutative Spreadability]
 {On non-commutative spreadability} 
\author{Maria Elena Griseta}
\address{Dipartimento di Matematica \\ 
Universit\`{a} degli Studi di Bari\\
Via E. Orabona, 4, 70125 Bari, Italy }
\email{mariaelena.griseta@uniba.it}



\date{May 26, 2023}

\begin{abstract}
We review some results on spreadable quantum stochastic processes and present the structure of some monoids acting on the index-set of all integers $\bz$. These semigroups are strictly related to spreadability, as the latter can be directly stated in terms of invariance with respect to their action.\\
We are mainly focused on spreadable, Boolean, monotone, and $q$-deformed processes. In particular, we give a suitable version of the Ryll-Nardzewski Theorem in the aforementioned cases.

\end{abstract}

\maketitle
\section*{Introduction}
Distributional symmetries for sequences of random variables have been intensely studied in classical probability, also for their essential applications to statistics and statistical mechanics \cite{Ka}. In recent years the development of non-commutative probability led to an intensive interest towards the study of distributional symmetries for the so-called quantum stochastic processes \cite{A1}.\\
These notes aim to review some recent results obtained for a prominent distributional symmetry, namely spreadability, and its relations with exchangeability and stationarity when applied to stochastic processes arising in some concrete models of non-commutative probability. 

\noindent A sequence $(X_n)_n$ of random variables is {\it spreadable} when for any $i_1,i_2, \ldots,i_n\in \bn$ the law of $(X_{i_1},X_{i_2},\ldots, X_{i_n})$ is the same of $(X_{g(i_1)},X_{g(i_2)},\ldots, X_{g(i_n)})$ for any strictly increasing map $g$ on the index set. The sequence is {\it exchangeable} if, instead, the same invariance property holds for all finite permutations $g$, and {\it stationary} when $g$ is exactly the one-side shift. This means that, in general, exchangeability implies spreadability, and this is stronger than stationarity.

\noindent We recall that a general version of de Finetti's theorem states that exchangeability for a sequence of random variables is equivalent to conditional independence and identical distribution with respect to the tail $\s$-algebra of the sequence itself \cite{Ka}. Moreover, the classical Ryll-Nardzewski Theorem asserts the equivalence between spreadability and exchangeability \cite{R}. This result is not always true in the non-commutative case. K\"{o}stler in \cite{K} showed examples of quantum stochastic processes in the so-called $W^*$-probability spaces which are spreadable while not being exchangeable. More recently, Crismale and Rossi in \cite{CriRosJFA} showed that in the CAR algebra there exist stationary states that are not spreadable. On the other hand, for some concrete models from non-commutative probability, such as Boolean, monotone, and $q$-deformed (with $|q| <1$) processes, spreadable and exchangeable processes still coincide. \cite{CFG2, CFG, CriRosAMPA}.

\noindent In these notes, we consider quantum stochastic processes indexed by a set $J$, generally identified with the set of integers $\bz$, on $C^*$-algebraic probability spaces. Here, assigning a quantum stochastic process having a $C^*$-algebra $\ga$ as sample space is equivalent to taking a state on the free product $C^*$-algebra $\ast_J \ga$ \cite{CrFid, CrFid2}. This means that if a distributional symmetry is implemented by an action of a group or a semigroup on $\ast_J \ga$, the invariance for the law of a process is transferred into the invariance of the corresponding state under the aforementioned action. 

\noindent The paper is organized as follows. In Section \ref{sec:prel} we give some definitions and notions concerning free product $C^*$-algebra and $C^*$-dynamical systems. In addition, we introduce some structures involved in spreadability, namely the monoid of strictly increasing maps on $\bz$, denoted by $\bl_\bz$, and its sub-monoid, $\bi_\bz$, generated by left and right hand-side partial shifts on the integers. We also recall the equivalence between unitarily equivalent classes of  $C^*$-stochastic processes on a sample algebra $\ga$ and states on the free product $C^*$-algebra of $\ast_J \ga$.  Furthermore, we present the definition of amenable semigroup, we introduce the monoid $\bj_\bz$ of strictly increasing maps on $\bz$ with cofinite range, and we see that it is properly included in the monoid $\bl_\bz$, whereas strictly contains $\bi_\bz$. We show that the spreadable states are exactly those invariant under the $*$-endomorphisms implemented by $\bj_\bz$. Finally, after denoting by $\bd_ \bz$ and $\be_ \bz$ the sub-monoids of $\bl_\bz$ generated by all forward and backward partial shifts respectively, we show that $\bj_\bz$ is semi-direct product of $\bd_\bz$ and $\be_\bz$, where the action of $\bz$ is generated by the conjugation by the one-step shift.\\
In Section \ref{sec:appl} we recall results given in \cite{CFG, CriRosAMPA, CFG2} concerning spreadability on the monotone, $q$-deformed and Boolean $C^*$-algebras, respectively. In particular, in Section \ref{sub:boolean} we pinpoint the Boolean case. Following \cite{CFL} and \cite{CrFid}, we present a suitable version of the Ryll-Nardzewski Theorem: the states invariant under the action of $\bj_\bz$ coincide with exchangeable states. In addition, one has that the latter are the same as the stationary states.

\section{Preliminaries}
\label{sec:prel}
For an arbitrary set $J$ and unital $C^*$-algebras $\{\ga_j\}_{j\in J}$, their unital free product $C^*$-algebra $\ast_{j\in J} \ga_j$ is the unique unital $C^*$-algebra, together with  unital monomorphisms $i_j:\ga_j\rightarrow \ast_{j\in J} \ga_j$ such that for any unital $C^*$-algebra $\gb$ and unital morphisms $\F_j:\ga_j\rightarrow \gb$, there exists a unique unital homomorphism $\F:\ast_{j\in J} \ga_j\rightarrow \gb$ making commutative the following diagram
\begin{equation*}
\xymatrix{ \ga_j \ar[r]^{i_j} \ar[d]_{\F_{j}} &
\ast_{j\in J} \ga_j \ar[dl]^{\quad\F\quad\quad \quad\quad} \\
\gb}
\end{equation*}
Here, we consider unital free product $C^*$-algebras based on a single unital $C^*$-algebra $\ga$, called {the free product
$C^*$-algebra}, and denoted by
$\ast_{J} \ga$. We refer the reader to \cite{CrFid, CrFid2, VDN} for further details.

Recall that the group of finite permutations $\bp_J$ on an arbitrary index-set $J$ is 
$$
\bp_J:=\bigcup\big\{\bp_\L\mid\L\,\text{finite subset of}\,J\big\},
$$
where $\bp_\L$ is the symmetric group associated with the finite set $\L$. If $J=\bz$, one can also consider the group generated by the one-step shift $\t(i):=i+1$
on $\bz$, which is canonically identified with $\bz$ itself. \\
If
$$
\bl_\bz:=\big\{f:\bz\to\bz\mid k<l\Rightarrow f(k)<f(l)\big\}\,,
$$
it is easy to check that it is a monoid with respect to composition. Recall that the $h$-{right hand-side partial shift}, $h\in\bz$, is the one-to-one map $\theta_h:\bz\rightarrow\bz$ such that
$$
\theta_h(k):=\left\{\begin{array}{ll}
k & \text{if}\,\, k<h\,, \\
k+1 & \text{if}\,\, k\geq h\,.
\end{array}
\right.
$$
Analogously, the $h$-{left hand-side partial shift}, $h\in\bz$, is the one-to-one map $\psi_h:\bz\rightarrow \bz$ such that
$$
\psi_h(k):=\left\{\begin{array}{ll}
k & \text{if}\,\, k>h\,, \\
k-1 & \text{if}\,\, k\leq h\,.
\end{array}
\right.
$$
We note that $\{\theta_h(k),\psi_h(k)\mid k\in\bz\}\subseteq\bl_\bz$.\\
Let us denote by $\bi_ \bz$ the sub-monoid of $\bl_\bz$ generated by all forward and backward partial shifts
$\{\th_h\}_{h\in\bz}$, and $\{\psi_h\}_{h\in\bz}$. We report, without proof, the following useful results (cf. \cite[Lemma 2.2]{CrFid}, and \cite[Proposition 2.1]{CFG}).
\begin{prop}
\label{spexi}
The following holds for a finite interval $[k,l]\subset\bz$.
\begin{itemize}
\item[(i)] There exists a cycle $\s_{k,l}\in\bp_\bz$ such that $\t([k,l])=\s_{k,l}([k,l])$.
\item[(ii)] For each $f\in\bl_\bz$ there exists $r_{k,l;f}\in\bi_\bz$ such that
$f([k,l])=r_{k,l;f}([k,l])$.
\end{itemize}
\end{prop}

The triplet $(\ga, M,\Gamma)$ is said to be a $C^*$-dynamical system if $\ga$ is a unital $C^*$-algebra, $M$ is a monoid, and  $\G: M\rightarrow {\rm End}(\ga)$ is a unital $*$-morphism. In some cases, $M$ is replaced by a group $G$, and in the $C^*$-dynamical system $(\ga, G,\a)$, $\a_g$ is a $*$-automorphism for each $g\in G$. In the latter case, one speaks of reversible dynamics, whereas dissipative dynamics appears in the absence of bijections, see, e.g. \cite{BR1}.

\noindent By $\cs(\ga)$, we denote the convex of the states on $\ga$, that is, the positive normalized linear functionals on $\ga$. $\cs(\ga)$ is weakly $*$-compact as $\ga$ is unital. The convex and compact in the $*$-weak topology, subset of all invariant states is
$$
\cs_M(\ga):=\big\{\f\in\cs(\ga)\mid\f\circ\G_g=\f\,,\,\,g\in M\big\}\,.
$$
The set of the extremal invariant states is denoted by $\ce_M(\ga)$, and it collects the {ergodic states} under the action $\G$ of $M$.

\noindent By universality, the groups $\bz$ and $\bp_\bz$ act in a natural way as $*$-automorphisms on the free product $C^*$-algebra $\ast_{\bz} \ga$ by shifting and permuting the indices of the generators, respectively. Moreover, it is possible to see (cf. \cite[Section 4]{CFG}) that there is an action by $*$-endomorphisms of both the monoids $\bi_\bz$ and $\bl_\bz$ on $\ast_{\bz} \ga$.\\
A state on $\ga$ is called symmetric if it belongs to $\cs_{\bp_\bz}(\ga)$, and shift invariant if it belongs to $\cs_{\bz}(\ga)$.\\
We have the following immediate consequences of Proposition \ref{spexi}.
\begin{cor}
\label{spexi1}
The following assertions hold:
\begin{itemize}
\item[(i)] $\cs_{\bp_\bz}(\ast_{\bz} \ga)\subseteq\cs_{\bz}(\ast_{\bz} \ga)$.
\item[(ii)] $\cs_{\bl_\bz}(\ast_{\bz} \ga)=\cs_{\bi_\bz}(\ast_{\bz} \ga)$.
\end{itemize}
\end{cor}
\noindent A quantum stochastic process with index-set $J$ is a quadruple $\big(\ga,\ch,\{\iota_j\}_{j\in J},\xi\big)$, where $\ga$ is a $C^{*}$-algebra, $\ch$ is an Hilbert space, the maps $\iota_j$ are $*$-homomorphisms of $\ga$ in $\cb(\ch)$, and
$\xi\in\ch$ is a unit vector, cyclic for the von Neumann algebra $\bigvee_{j\in J}\iota_j(\ga)$.
\begin{defn}
For the $C^*$-algebra $\ga$, $n\in\bn$, $j_1,\ldots j_n\in\bz$, $A_1,\ldots A_n\in\ga$,
the stochastic process $\big(\ga,\ch,\{\iota_j\}_{j\in \bz},\xi\big)$ is said to be
\begin{itemize}
\item[{\bf-}] {stationary} if
$$
\langle\iota_{j_1}(A_1)\cdots\iota_{j_n}(A_n)\xi,\xi\rangle
=\langle\iota_{j_{1}+1}(A_1)\cdots\iota_{j_{n}+1}(A_n)\xi,\xi\rangle\,;
$$
\item[{\bf-}] {exchangeable} if for each $g\in\bp_{\bz}$,
\begin{equation*}
\langle\iota_{j_1}(A_1)\cdots\iota_{j_n}(A_n)\xi,\xi\rangle
=\langle\iota_{g(j_{1})}(A_1)\cdots\iota_{g(j_{n})}(A_n)\xi,\xi\rangle\,;
\end{equation*}
\item[{\bf-}] {spreadable}
if for each $g\in\bl_\bz$,
\begin{equation}
\label{spre}
\langle\iota_{j_1}(A_1)\cdots\iota_{j_n}(A_n)\xi,\xi\rangle
=\langle\iota_{g(j_{1})}(A_1)\cdots\iota_{g(j_{n})}(A_n)\xi,\xi\rangle\,.
\end{equation}
\end{itemize}
\end{defn}
Condition \eqref{spre} and Proposition \ref{spexi} entail that spreadable stochastic processes are exactly those whose joint distributions are invariant under all partial shifts. By \cite[Theorem 3.4]{CrFid}, the states on $\ast_{\bz} \ga$ uniquely correspond to quantum stochastic processes on the sample algebra $\ga$. More in detail, for the quadruple $\big(\ga,\ch,\{\iota_j\}_{j\in \bz},\xi\big)$ one finds a state $\f\in\cs\big(\ast_{\bz}\ga\big)$, such that 
$$
\f(i_{j_1}(a_1)i_{j_2}(a_2)\cdots i_{j_n}(a_n))=\langle \iota_{j_1}(a_1)\iota_{j_2}(a_2)\cdots\iota_{j_n}(a_n)\xi,\xi\rangle
$$
for every $n\in\bn$, $j_1,j_2,\ldots, j_n\in\bz$, and $a_1,a_2,\ldots,a_n\in\ga$. Conversely, for $\f\in\cs\big(\ast_{\bz}\ga\big)$, if $(\pi_\f,\ch_\f,\xi_\f)$ denotes its GNS representation, one recovers a stochastic process after defining $\iota_j(a):=\pi_\f(i_j(a))$, for all $j\in\bz$ and $a\in\ga$. Therefore,
\begin{itemize}
    \item[{\bf-}] Stationary processes correspond to shift invariant states on $*_\bz\ga$;
    \item[{\bf-}] Exchangeable processes correspond to symmetric states on $*_\bz\ga$;
    \item[{\bf-}] Spreadable processes correspond to spreading invariant states on $*_\bz\ga$.
\end{itemize}
Note that Corollary \ref{spexi1} yields the following information:
\begin{itemize}
\item[{\bf-}] If a process is exchangeable, then it is stationary;
\item[{\bf-}] A process is spreadable if and only if it is invariant under the action of  $\bi_\bz$.
\end{itemize}
In the following sections, we will show that, in particular cases, spreadability can be studied by a suitable monoid included between $\bl_\bz$ and $\bi_\bz$.\\
Here we recall the definition of amenability for semigroups.\\
A discrete semigroup $S$ is said to be left (or right) amenable if there exists a state $\varphi$ on $\ell^{\infty}(S)$ such that $\varphi(l_sf) =\varphi(f)$ (or $\varphi(r_sf) =\varphi(f)$), for every $s \in S$ and $f\in\ell^{\infty}(S)$, where $l_sf(t):=f(st)$ (or $r_sf(t):=f(ts)$), for any $t \in S$. We observe that left amenability and right amenability are not the same notion. \\
In addition, let us denote by $\bd_ \bz$ and $\be_ \bz$ the sub-monoids of $\bl_\bz$ generated by all forward and backward partial shifts
$\{\th_h\}_{h\in\bz}$ and $\{\psi_h\}_{h\in\bz}$, respectively. In addition, let us take
$$
\bj_\bz:=\big\{f\in\bl_\bz : |\bz\smallsetminus f(\bz)|<+\infty\big\}\,.
$$
$\bj_ \bz$ is also a sub-monoid of $\bl_\bz$. If $f\in\bl_\bz$, we denote $\d_f:=\bz\smallsetminus f(\bz)$, and note that
\begin{itemize}
\item[{\bf -}] $f\in\bj_\bz\iff\d_f<+\infty$,
\item[{\bf -}] $f=\t^n$ for some $n\in\bz\iff\d_f=\emptyset$.
\end{itemize}
By using Lemma 1 of \cite{CFG2} one achieves the following properties:
$$\cs_{\bl_\bz}(\ast_{\bz} \ga)=\cs_{\bj_\bz}(\ast_{\bz} \ga)=\cs_{\bi_\bz}(\ast_{\bz} \ga)\,,
$$
$$
\bi_ \bz\subsetneq\bj_ \bz\subsetneq\bl_ \bz\,.
$$
The fact that spreadability can be expressed in terms of invariance with respect to the action of $\bj_ \bz$ will be useful in the case of Boolean processes, as we will see in Section \ref{sub:boolean}. \\
finally, recall that for two unital semigroups $M$ and $N$,  suppose that $M$ acts on $N$ by the morphisms
$$
M\ni m\mapsto\eta_m\in{\rm Mor}(N)\,.
$$
The semi-direct product $M\,_{\eta}\!\!\ltimes N$ is defined as follows.
As a set, $M\,_{\eta}\!\!\ltimes N:=M\times N$, whereas a binary operation is given by
\begin{equation}
\label{mno}
(m_1,n_1)(m_2,n_2):=\big(m_1m_2, n_1\eta_{m_1}(n_2)\big)\,.
\end{equation}
It is easy to check that $M\times N$, equipped with the multiplicative law in Equation \eqref{mno} defines a monoid whose unit is $e_{M\,\!_{\eta}\!\ltimes N}=(e_M,e_N)$.\\
The following result, which is Proposition 4 in \cite{CFG2}, allows us to see $\bj_\bz$ as a semi-direct product of semigroups.
\begin{prop}
\label{propmonoid2}
If for each $m\in\bz$ one defines $\eta_m(\cdot):= \t^m\cdot \t^{-m}$, it follows
$$
\bj_\bz=\bz\,_{\eta}\!\!\ltimes\bd_\bz=\bz\,_{\eta}\!\!\ltimes\be_\bz\,.
$$
\end{prop}

\section{Applications to the non-commutative setting}
\label{sec:appl}
In this section we recall several results about spreadability in the non-commutative probability field. In particular, we deal with spreadability for monotone and Boolean commutation relations, for the $C^*$-algebra generated by the annihilators acting on the $q$-deformed Fock space, and finally for the {\rm CAR} algebra. 

\subsection{Spreadability for Monotone Commutation Relations}
\label{sub:mono}
In this subsection, after recalling several notions concerning the monotone Fock space, we recall that the monoids generated by partial shifts on $\bz$ act as $*$-endomorphism on the concrete unital $C^*$-algebra generated by the annihilators acting on the monotone Fock space. Finally, we report the characterization of spreading invariant states with respect to such an action. We refer the reader to \cite{CFG} and the references cited therein for a detailed treatment.

For $k\geq 1$, denote $I_k:=\{(i_1,i_2,\ldots,i_k) \mid i_1< i_2 < \cdots <i_k, i_j\in \mathbb{Z}\}$. The discrete monotone Fock space is the Hilbert space $\cf_m:=\bigoplus_{k=0}^{\infty} \ch_k$, where for any $k\geq 1$, $\ch_k:=\ell^2(I_k)$, and $\ch_0=\mathbb{C}\Om$, $\Om$ being the Fock vacuum. Each
$\ch_k$ is called the $k$th-particle space and we denote by $\cf^o_m$ the total set of finite particle vectors in $\cf_m$. \\
Let $(i_1,i_2,\ldots,i_k)$ be an increasing sequence of integers. The generic element of the canonical basis of $\cf_m$ is denoted by $e_{(i_1,i_2,\ldots,i_k)}$, with the convention $e_\emptyset=\Om$. We often write $e_{(i)}=e_i$ to ease the notations.
The monotone creation and annihilation operators are respectively given, for any $i\in \mathbb{Z}$, by
\begin{equation*}
a^\dagger_i e_{(i_1,i_2,\ldots,i_k)}:=\left\{
\begin{array}{ll}
e_{(i,i_1,i_2,\ldots,i_k)} & \text{if}\, i< i_1\,, \\
0 & \text{otherwise}\,, \\
\end{array}
\right.
\end{equation*}
\begin{equation*}
a_ie_{(i_1,i_2,\ldots,i_k)}:=\left\{
\begin{array}{ll}
e_{(i_2,\ldots,i_k)} & \text{if}\, k\geq 1\,\,\,\,\,\, \text{and}\,\,\,\,\,\, i=i_1\,,\\
0 & \text{otherwise}\,. \\
\end{array}
\right.
\end{equation*}
One can check that both $a^\dagger_i$ and $a_i$ have unital norm (\cite[Proposition 8]{Boz}), are mutually adjoint and satisfy the following relations
\begin{equation*}
\begin{array}{ll}
  a^\dagger_ia^\dagger_j=a_ja_i=0 & \text{if}\,\, i\geq j\,, \\
  a_ia^\dagger_j=0 & \text{if}\,\, i\neq j\,.
\end{array}
\end{equation*}
In addition, the following commutation relation
\begin{equation*}
a_ia^\dagger_i=I-\sum_{k\leq i}a^\dagger_k a_k
\end{equation*}
is also satisfied, where the convergence of the sum is in the strong operator topology, see \cite[Proposition 3.2]{CFL}. \\
We denote by $\gam$ and $\gam_o$ the concrete unital $C^*$-algebra, together with its dense unital $*$-algebra generated by the annihilators $\{a_i\mid i\in\mathbb{Z}\}$, acting on the monotone Fock space. Both algebras are canonically $*$-subalgebras of $\cb(\cf_m)$.\\
We recall the structure of a Hamel basis for $\gam_o$, useful to get results on spreadability. We denote by $\L$ the index set such that $\{X_\l\}_{\l\in\L}$ denotes the  $\l$-forms, that is words of type
\begin{equation*}
X_\l=a_{i_1^{(\l)}}^{\dagger}\cdots a_{i_{m(\l)}^{(\l)}}^{\dagger} a_{j_1^{(\l)}}\cdots a_{j_{n(\l)}^{(\l)}}\,,
\end{equation*}
for $i_1^{(\l)}<i_2^{(\l)}<\cdots < i_{m(\l)}^{(\l)}, j_1^{(\l)}>j_2^{(\l)}> \cdots > j_{n(\l)}^{(\l)}$, $m(\l),n(\l)\geq 0$, where, as usual, $m(\l)=n(\l)=0$ corresponds to the identity $I$. We refer to $m(\l)+n(\l)$ as the length of the word. As  all the $\l$-forms like $a^\dag_ia_i$ are in one-to-one correspondence with $\bz$, we denote the complement set of them by $\G$, namely $\G:=\L\setminus\bz$.\\
The families $(X_\l)_{\l\in\G}$ and $(a_ia^\dag_i)_{i\in \mathbb{Z}}$ make up a Hamel basis for the monotone $*$-algebra (see \cite[Theorem 3.4]{CFG}).\\ 
One sees that $\bi_\bz$ act, by $*$-endomorphisms on the whole $\gam$. This gives that the spreading invariant states allow to determine spreading invariant processes on the monotone $C^*$-algebra.   \\
\noindent To this aim, after denoting by $\ga_o$ the linear span of words in $\gam_o$ with positive length, and by $\ga$ its norm closure, from Corollary 5.10 of \cite{CFL}, we recall that $\gam$ results to be the $C^*$-algebra obtained by adding the identity to $\ga$. Then any $Y\in\gam$ decomposes as $Y:=X+cI$, where $X\in\ga$ and $c\in\mathbb{C}$.\\
On $\gam$ the state at infinity $\om_\infty$ is well defined as 
$$
\om_\infty(Y)=\om_\infty(X+cI):=c\,,
$$
and the vacuum state $\om$ is given as usual as
$$
\om(Y):=\langle Y \Om,\Om\rangle\,.
$$
The following result yields the simplex of spreading-invariant states in the monotone case (see \cite[Proposition 4.5]{CFG}).
\begin{prop}
\label{inv}
The $*$-weakly compact set of spreading invariant states on $\gam$ is
$$
\cs_{\mathbb{I}_\mathbb{Z}}(\gam)=\{(1-x)\om_\infty + x\om \mid x\in [0,1]\}\,.
$$
\end{prop}

\subsection{Spreadability for q-deformed stochastic processes}
\label{sub:qdeformed}
This subsection is devoted to recalling that the Ryll-Nardzewski Theorem holds on the concrete $C^*$-algebra generated by the annihilators acting on the $q$-deformed Fock space if $|q|<1$. In particular, the Fock vacuum vector state is the unique spreading-invariant, shift invariant, and stationary state.  We refer the reader to \cite{CriRosAMPA} for a more detailed discussion.

\noindent We start with the definition of $q$-deformed Fock space \cite{BS}. \\
Let $-1<q<1$ and $\ch$ a Hilbert space. The $q$–deformed Fock space $\G_q(\ch)$ is the completion of the algebraic linear span of the vacuum vector $\zeta_q$, together with vectors.
$$
f_1\otimes \cdots \otimes f_n\,, \quad f_j\in\ch\,, j=1,\ldots,n\,, n=1,2,\ldots
$$
with respect to the $q$–deformed inner product
$$
{\langle f_1 \otimes \cdots f_m,g_1 \otimes \cdots g_n\rangle}_q:=\d_{n,m}\sum_{\pi\in\bp_n}q^{i(\pi)}\langle f_1,g_{\pi(1)}\rangle_{\ch}\cdots\langle f_n,g_{\pi(n)}\rangle_{\ch}\,,
$$
where $i(\pi)$ is the number of inversions of $\pi\in\bp_n$, the symmetric group of order $n$. Fix $f,f_1,\ldots f_n\in\ch$. Define the creator $l^{\dag}(f)$ as
$$
l^{\dag}(f)\zeta_q=f\,, \quad l^{\dag}(f)f_1\otimes\cdots\otimes f_n=f\otimes f_1\cdots\otimes f_n\,,
$$
and the annihilator $l(f)$ as
\begin{equation*}
\begin{split}
   l(f)\zeta_q &=0\,,\\
   l(f)f_1\otimes\cdots\otimes f_n &= \sum_{k=1}^n q^{k-1}\langle f_k,f\rangle_{\ch} f_1 \otimes\cdots f_{k-1}\otimes f_{k+1}\otimes\cdots\otimes f_n\,.
\end{split}
\end{equation*}
 $l^{\dag}(f)$ and $l(f)$ are extended by linearity on the dense subspace of finite particle vectors of $\G_q(\ch)$, where they are also continuous. So their extensions are in $\mathcal{B}(\G_q(\ch))$, are mutually adjoint with respect to the $q$–deformed inner product, and they are the Fock representation of the $q$-commutation relations, namely
$$
l(f)l^{\dag}(g)-ql^{\dag}(g)l(f)=\langle g,f\rangle_{\ch}I\,, \quad f,g\in\ch\,,
$$
$I$ being the identity on $\mathcal{B}(\G_q(\ch))$.\\
Let us now take $\ch =\ell^2(\bz)$ with the canonical orthonormal basis $\{e_j\}_{j\in\bz}$, and on $\G_q(\ell^2(\bz))$ denote $l_j := l(e_j)$, $j \in\bz$. The concrete $C^*$–algebra $\mathfrak{R}_q$ and its subalgebra $\mathfrak{S}_q$, acting on $\G_q(\ell^2(\bz))$ are the unital $C^*$–algebras generated by the annihilators $\{l_j | j \in \bz\}$, and the self-adjoint part of annihilators $s_j := l_j + l^{\dag}_j$, $j\in\bz$, respectively. If $\omega_q :=\langle \cdot \z_q,\z_q\rangle$ denotes the Fock vacuum vector state, one has the following
\begin{thm}
\label{spreadq}
$$
\cs_{\bp_\bz}(\mathfrak{R}_q)=\cs_{\bl_\bz}(\mathfrak{R}_q)=\cs_{\bz}(\mathfrak{R}_q)=\{\omega_q\}\,.
$$
and
$$
\cs_{\bp_\bz}(\mathfrak{S}_q)=\cs_{\bl_\bz}(\mathfrak{S}_q)=\cs_{\bz}(\mathfrak{S}_q)=\{\omega_q\}\,.
$$
hold.
\end{thm}
We omit the proof. The reader is referred to \cite[Theorem 3.3 and Corollary 3.4]{DF} and \cite[Proposition 6.2]{CrFid} for more details.

\subsection{Spreading Invariant States on the Boolean Algebra}
\label{sub:boolean}
We report, without proof, the main results obtained in \cite{CFG2}, in which the authors proved in a direct way that the monoid $\bj_\bz$ acts on the concrete Boolean unital $C^*$-algebra by unital $*$-endomorphisms. In the last part of the section, we give a suitable Ryll-Nardzewski Theorem in this setting. 

\noindent We start by recalling the definition of Boolean Fock space.\\
Let $\ch$ be a complex Hilbert space. The Boolean Fock space $\cf_{\rm boole}(\ch)$ over $\ch$ is the direct sum $\cf_{\rm boole}(\ch)=\bc\oplus\ch$, and the vacuum vector is
$e_{\#}:=1\oplus 0$. The vacuum vector state is denoted by $\om_{\#}:=\langle{\bf\cdot}\,e_{\#},e_{\#}\rangle$.\\
For $\g\in\bc$ and $f,g\in\ch$, the creation and annihilation operators are defined as follows:
\begin{equation*}
b^{\dag}(f)(\g\oplus g):= 0\oplus \g f, \,\,\,\, b(f)(\g\oplus g):= \langle g,f\rangle\oplus 0\,.
\end{equation*}
They are mutually adjoint and bounded.\\
The concrete Boolean $C^*$-algebra $\gpb_\ch$ is that generated by all creators and the identity $I$.
Since the $*$-algebra generated by the $b^{\dag}$ consists of all finite rank operators, we easily get
\begin{equation}
\label{decom}
\gpb_\ch=\ck(\bc\oplus\ch)+\bc I\,,
\end{equation}
where $\ck(\bc\oplus\ch)$ denotes the $C^*$-algebra of compact linear operators on $\ch$.
Here, we deal with the case $\ch=\ell^2(\bz)$, where the canonical basis is $(e_j)_{j\in\bz}$. Therefore,
\begin{align*}
\cf_{\rm boole}(\ell^2(\bz))=&\bc e_\#\oplus\ell^2(\bz)=\ell^2(\{\#\}\sqcup\bz)\,,\\
\gpb_{\ell^2(\bz)}=&\ck(\ell^2(\{\#\}\sqcup\bz))+\bc I\,.
\end{align*}
With the notations $b_j := b(e_j)$, $b^{\dag}_j := b^{\dag}(e_j)$, we can see that the following Boolean commutation relation holds for arbitrary $i,j\in\bz$ 
\begin{equation*}
b_ib^\dag_j   = \d_{i,j}\bigg(I -\sum_{k\in\bz}b^\dag_k b_k\bigg)\,.
\end{equation*}
From now on, we use the shorthand notation $\bz_\#:=\{\#\}\sqcup\bz$.\\
For each subset $A\subset\bz_\#$, $P_A\in\cb(\ell^2(\bz_\#))$ will denote the self-adjoint projection onto the closed subspace of $\ell^2(\bz_\#)$ generated by the $e_j$, $j\in A$. For the canonical system of matrix--units $\{\varepsilon_{ij}\mid i,j\in \bz_{\#}\}$ in $\cb(\ell^2(\bz_\#))$, one has
$$
b_j=\varepsilon_{\# j}\,,\,\, b^\dagger_j=\varepsilon_{j\#}\,,\,\,
b_ib^\dagger_j=\varepsilon_{\#\#}\d_{i,j}\,,\,\,b^\dagger_ib_j=\varepsilon_{ij}\,,\quad i,j\in\bz\,.
$$
It is well known that the groups $\bz$ and $\bp_\bz$ naturally act on $\gpb_{\ell^2(\bz)}$ by $*$-automorphisms ({e.g.}, \cite{CrFid, CrFid2, CFL}).\\
For the action of $\bj_{\bz}$ on $\gpb_{\ell^2(\bz)}$, we consider $\mathcal{K}_o(\ell^2(\bz_\#))$, the $*$-algebra of finite rank operators on the Boolean Fock space. On the dense $*$-algebra $\mathcal{K}_o(\ell^2(\bz_\#))+\bc I$, dense in the norm topology in $\gpb_{\ell^2(\bz)}$, we define
\begin{equation*}
\a^{(o)}: f\in\bj_ \bz\mapsto \a^{(o)}_{f}\in\big(\mathcal{K}_o(\ell^2(\bz_\#))+\bc I)^{\mathcal{K}_o(\ell^2(\bz_\#))+\bc I}
\end{equation*}
such that
\begin{equation}
\label{alfa}
\begin{split}
&\a^{(o)}_f(I):=I\,,\\
&\a^{(o)}_f(\eps_{kl}) :=\eps_{f_\#(k)f_\#(l)}, \quad k,l\in\bz_\#\,, \\
\end{split}
\end{equation}
where
$$
f_\#(k):=\left\{
\begin{array}{ll}
f(k) & \text{if $k\in\bz$}\,, \\
k & \text{if $k=\#$}\,.
\end{array}
\right.
$$
The $\a^{(o)}_f$ are well defined because $\big\{\eps_{kl}\mid k,l\in\bz_\#\big\}\subset\mathcal{K}_o(\ell^2(\bz_\#))$ is a Hamel basis. In \cite[Theorem 1]{CFG2},
the authors proved that $\a^{(o)}$ extends to an action of $\bj_\bz$ by unital $*$-endomorphisms of the boolean $C^*$-algebra $\gpb_{\ell^2(\bz)}$. We briefly report the results here.

\noindent Indeed, for any $f\in\bj_{\bz}$ we define
\begin{equation*}
V_fe_k:=e_{f_\#(k)}\,, \quad k\in\bz_\#\,,
\end{equation*}
which extends to an isometry on the whole $\ell^2(\bz_\#)$, denoted again by $V_f$ (see \cite[Lemma 2]{CFG2}). By identity \eqref{decom}, any $X\in\gpb_{\ell^2(\bz)}$ is decomposed as $X:=K+\g I$, where $K\in\mathcal{K}(\ell^2(\bz_\#))$ and $\g\in\bc$. Thus, the state at infinity $\om_\infty$ is also in this setting well defined as
\begin{equation*}
\om_\infty(X):=\g\,.
\end{equation*}
For $f\in\bj_\bz$, define the linear maps
\begin{equation}
\label{act}
\a_f(X)=V_fXV_f^{*}+\om_{\infty}(X)P_{\d_f}\,,\quad X\in\mathfrak{b}_{\ell^2(\bz)}\,.
\end{equation}
\begin{thm}
\label{thmalpha}
The map
$$
\a:f\in\bj_\bz\rightarrow\a_f\in{\rm End}\big(\mathfrak{b}_{\ell^2(\bz)}\big)
$$
given in Equation \eqref{act} provides a representation of the monoid $\bj_\bz$ in ${\rm {\rm End}}\big(\mathfrak{b}_{\ell^2(\bz)}\big)$ extending the linear maps given in Equation \eqref{alfa}.
\end{thm}

The following result establishes an equivalence between stationary, spreadable and exchangeable processes on the concrete Boolean $C^*$-algebra, thus, in particular, realising a version of Ryll-Nardzewski Theorem \cite{R} in the Boolean setting. We refer the reader to \cite[Proposition 5]{CFG2} for a proof.
\begin{prop}
\label{spr}
The states on $\mathfrak{b}_{\ell^2(\bz)}$ which are spreading invariant, coincide with the stationary and symmetric states:
\begin{align*}
\cs_{\bj_\bz}\big(\mathfrak{b}_{\ell^2(\bz)}\big)=&\cs_{\bp_\bz}\big(\mathfrak{b}_{\ell^2(\bz)}\big)
=\cs_\bz\big(\mathfrak{b}_{\ell^2(\bz)}\big)\\
=&\{\l\om_{\#} +(1-\l)\om_{\infty} \mid \l\in[0,1] \}\,.
\end{align*}
\end{prop}

\subsection{Remarks on the CAR algebra}
\label{sub:CAR}
The Canonical Anticommutation Relations ({\rm CAR} for short) algebra over $\bz$ is the universal unital $C^*$-algebra {\rm CAR}($\bz$), with unit $I$, generated by the set $\{A_j, A_j^{\dag} \mid j\in\bz\}$ (i.e. the Fermi annihilators and creators respectively), satisfying the relations
\begin{align*}
    (A_j)^{*} = A_j^{\dag}\,, & \quad \{A_j^{\dag},A_k\}=\delta_{j,k}I\,, & \{A_j,A_k\}=\{A_j^{\dag},A_k^{\dag}\}=0\,,\quad j,k\in\bz\,,
\end{align*}
where $\{\cdot,\cdot\}$ is the anticommutator and $\delta_{j,k}$ is the Kronecker symbol.\\
Note that by definition
$$
{\rm CAR}(\bz) =\overline{{\rm CAR}_0(\bz)} \,,
$$
where
$$
{\rm CAR}_0(\bz) := \bigcup\{{\rm CAR}(F) : F\subset \bz \, \text{finite}\}
$$
is the (dense) subalgebra of the localized elements, and ${\rm CAR}(F)$ is the $C^*$-subalgebra generated by the finite set $\{A_j,A_j^{\dag}:j\in F\}$.
\begin{thm}
\label{CAR}
The following inclusion holds
$$
\cs_{\bj_\bz}({\rm CAR}(\bz))\subsetneq \cs_{\bz}({\rm CAR}(\bz))\,.
$$
\end{thm}
\noindent In \cite[Proposition 4.1., Theorem 4.2.]{CriRosJFA} the authors found a state $\omega$ on {\rm CAR}($\bz$) such that
\begin{equation*}
\omega(A^{\dag}_m A_n)=i\frac{3C}{\pi^2(m-n)^2}\,,
\end{equation*}
where $i$ denote the imaginary unity of $\bc$, $m,n\in\bz$ with $m>n$ and $C$ is a positive constant. One figures out that the state above is stationary but not spreadable. \\
Finally, we mention the case concerning the self-adjoint part of ${\rm CAR}(\bz)$. This is by definition the unital $C^*$-algebra generated by the position operators, say $\mathfrak{C} :=C^{*}(x_j:j\in\bz)$, where $x_j:=A_j+A_j^{\dag}$ for every $j\in\bz$. The $x_j$'s anticommute with one another and their square is the identity, that is
\begin{center}
$x_jx_k + x_kx_j = 0$, for all $j\neq k$, and  $x_j^2=I$, for all $j\in\bz$.
\end{center}
On $\mathfrak{C}$ as well, one finds that there exist stationary states that are not spreadable.

\subsection*{Acknowledgment}
The author is partially supported by Progetto GNAMPA 2023 ``Metodi di Algebre di Operatori in Probabilità non commutativa'' and by Italian PNRR MUR project PE0000023-NQSTI, CUP H93C22000670006.

\end{document}